\documentclass[11pt,letterpaper]{amsart}
\usepackage{graphicx}
\usepackage{amsmath,amssymb,amsfonts,mathrsfs}
\usepackage{amsthm,enumerate,multirow,multicol}

\usepackage{psfrag,color}
\usepackage[small,hang,bf,up,sf]{caption}

\newcounter{lletres}

\newtheorem{Talpha}[lletres]{Theorem}
\newtheorem{Lalpha}[lletres]{Lemma}

\newtheorem{thm}{Theorem}
\newtheorem{defi}[thm]{Definition}

\newtheorem{lma}[thm]{Lemma}

\newcommand{\D}{\mathbb D}

\newcommand{\N}{\mathbb N}

\newcommand{\U}{\mathcal U}
\newcommand{\B}{\mathcal B}

\newcommand{\Hi}{H^\infty}
\newcommand{\Hp}{H^p}
\newcommand{\eps}{\varepsilon}
\newcommand{\tdel}{\widetilde\Delta}
\newcommand{\hdel}{\widehat{\Delta}}

\begin{document}

\title[Sampling and interpolation for analytic selfmappings...]{Sampling and interpolation for analytic selfmappings of the disc.}
\author[Nacho Monreal Gal\'{a}n]{Nacho Monreal Gal\'{a}n}
\address{Nacho Monreal Gal\'{a}n, Departament of Mathematics, University of Crete, Voutes Campus, 70013 Heraklion, Crete, Greece.}
\email{nacho.mgalan@gmail.com}
\author[Michael Papadimitrakis]{Michael Papadimitrakis}
\address{Michael Papadimitrakis, Departament of Mathematics, University of Crete, Voutes Campus, 70013 Heraklion, Crete, Greece.}
\email{mihalis.papadimitrakis@gmail.com}

\thanks{
The first author is supported in part by the projects MTM2011-24606 and by the research project PE1(3378) implemented within the framework of the Action 
``Supporting Postdoctoral Researchers'' of the Operational Program ``Education and Lifelong Learning'' (Action's Beneficiary: General Secretariat for 
Research and Technology), co-financed by the European Social Fund (ESF) and the Greek State.}

\begin{abstract}
Two different problems are considered here. First, a characterization of sampling sequences for the class of analytic functions from the disc into itself is given. Second, a version of Schwarz-Pick Lemma for $n$ points leads to an interpolation problem for the same functions, which may be considered as a particular case of the classical Nevanlinna-Pick interpolation problem.
\end{abstract}

\maketitle

\section{Notation and some basic facts.}
 
Let $\Hi$ be the Banach space of bounded analytic functions in the open unit disc $\D$ of the complex plane, with norm
$$\|f\|_\infty = \sup \{|f(z)| : z \in \D\}<\infty.$$
Let also $\U$ be the set of functions $f \in \Hi$ with $\|f\|_\infty \leq 1$. Thus, $\U$ is just the closed unit ball of $\Hi$.

The pseudohyperbolic distance between $z,w\in\D$ is 
$$\rho(z,w)=\Big|\frac{z-w}{1-\overline w z}\Big|.$$
Then the hyperbolic distance between $z$ and $w$ in $\D$ is 
$$\beta(z,w)=\log\frac{1+\rho(z,w)}{1-\rho(z,w)},$$
and the hyperbolic disc with center $z_0$ and radius $R>0$ is given by
$$D_h(z_0,R)=\{z\in\D : \beta(z,z_0)<R\}.$$
The corresponding hyperbolic circle is denoted $C_h(z_0,R)$.

The Schwarz-Pick Lemma asserts that, if $f\in\U$ and $z,w\in\D$, then 
$$\beta(f(z),f(w))\leq\beta(z,w).$$ 
Furthermore, this inequality becomes equality for a pair (and then for all pairs) of distinct points of the disc if and only if $f$ is an automorphism of the disc, i.e. if $f(z)=\lambda\frac{z-z_0}{1-\overline{z_0} z}$ for some $z_0\in\D$ and some $\lambda$ with $|\lambda|=1$ (see for instance \cite{Gar}). 

\section{Sampling sequences for $\U$.}

In \cite{BN} the authors defined the concept of sampling sequence for the Bloch space and gave a characterization of it. We say that an analytic function $f$ in $\D$ is in the Bloch space $\B$ if
$$\|f\|_{\B}=\sup_{z,w\in\D}\frac{|f(z)-f(w)|}{\beta(z,w)}<\infty.$$
The quantity $\|f\|_{\B}$ defines a semi-norm and $\B$ becomes a Banach space with the norm $\|f\|=|f(0)|+\|f\|_{\B}$. It is well known that
$$\|f\|_{\B}=\sup_{z\in\D}(1-|z|^2)|f'(z)|.$$

A sequence $Z=\{z_n\}$ of distinct points in $\D$ is said to be sampling for $\B$ if there exists $\delta>0$ such that for every $f\in\B$ one has that
$$\sup_{n\neq m}\frac{|f(z_n)-f(z_m)|}{\beta(z_n,z_m)}\geq \delta\|f\|_{\B}.$$

In the cited paper it is proved that $Z$ is a sampling sequence if and only if $Z$ is $R$-dense in $\D$ for some $R>0$, that is, $Z\cap D_h(z,R)\neq\emptyset$ for every $z\in\D$. In \cite{Seip} one can find a deeper description of this problem.

For an analytic function $f$ in $\U$ we define
$$N(f)=\sup_{z,w\in\D}\frac{\beta(f(z),f(w))}{\beta(z,w)}$$
and the Schwarz-Pick Lemma says that $N(f)\leq 1$. Again, it is well known that
$$N(f)=\sup_{z\in\D}\frac{(1-|z|^2)f'(z)}{1-|f(z)|^2}.$$

Focusing on a similarity between the space $\B$ and the class of functions $\U$, we see that in both of them the functions satisfy a Lipschitz type inequality. In the case of a Bloch function $f$, the inequality involves Euclidean distance in the image space and the Lipschitz constant is $\|f\|_{\B}$. In the case of an analytic self-mapping $f$ of the disc, the inequality involves hyperbolic distance in the image space and the Lipschitz constant is $N(f)$. Therefore, a sampling problem for $\U$ may be understood as the hyperbolic version of a sampling problem for $\B$.

\begin{defi}
A sequence $Z=\{z_n\}$ of distinct points in $\D$ is said to be sampling for $\U$ if there exists $\delta>0$ such that for every $f\in\U$ one has that
$$\sup_{n\neq m}\frac{\beta(f(z_n),f(z_m))}{\beta(z_n,z_m)}\geq \delta N(f).$$
\end{defi}

As in the case of $\B$, this definition is conformally invariant, in the sense that if $Z=\{z_n\}$ is a sampling sequence for $\U$ and $\tau$ is an automorphism of the disc then the sequence $\tau(Z)=\{\tau(z_n)\}$ is also a sampling sequence with the same sampling constant $\delta$.

We shall prove the following characterization of sampling sequences for $\U$.

\begin{thm}
A sequence $Z$ of distinct points in $\D$ is a sampling sequence for $\U$ if and only if $Z$ is $R$-dense in $\D$ for some $R>0$. 
\end{thm}

This characterization coincides with the one which holds for $\B$. The proof of the sufficiency though is different from the one  given in \cite{BN}.

\begin{proof}
	
\noindent \textit{Necessity.}
	
\noindent Let $Z=\{z_n\}$ be a sampling sequence for $\U$. Then there is $\delta>0$ such that for every $f\in\U$ and hence for $f(z)=\frac z2$, 
$$\sup_{n\neq m}\frac{\beta(f(z_n),f(z_m))}{\beta(z_n,z_m)}\geq\delta N(f).$$
Now, it is easy to see that 
$$\frac{\beta\big(\frac z2,\frac w2\big)}{\beta(z,w)}\to 0\quad \text{if}\quad |z|\to 1\,\,\text{or}\,\,|w|\to 1.$$ 
So there exists $r>0$ which depends only on $C$ so that $|z_n|<r$ for some $n$. Therefore, if $R=\log\frac{1+r}{1-r}$, then the hyperbolic disc $D_h(0,R)$ intersects $Z$. By conformal invariance of a sampling sequence, this implies that $Z$ is $R$-dense. 
	
\noindent \textit{Sufficiency.}
	
\noindent Now assume that the sequence $Z=\{z_n\}$ of distinct points in $\D$ is $R$-dense. Arguing by contradiction, we suppose that $Z$ is not a sampling sequence for $\U$. Hence, for every $k\in\N$ there exists $f_k\in\U$ such that
$$\sup_{n\neq m}\frac{\beta(f_k(z_n),f_k(z_m))}{\beta(z_n,z_m)}<\frac{1}{k}\,N(f_k).$$
Using appropriate automorphisms of $\D$, we may assume that 
$$f_k(0)=0,\quad |f_k'(0)|\geq\frac 12\,N(f_k).$$
Now we remove all points of $Z$ which belong to $D_h(z_1,R)$ except $z_1$. Of the points of $Z$ which are outside of $D_h(z_1,R)$ we take the one with the smallest index, say $z_m$, and we remove all points of $Z$ which belong to $D_h(z_m,R)$ except $z_m$. We continue inductively producing a new sequence $Z'$ smaller than $Z$. It is clear that the hyperbolic discs with centers from $Z'$ and radius $\frac R2$ are disjoint. Now let us take any $z\in\D$. Then $D_h(z;R)$ contains some $z_k\in Z$. If $z_k\in Z'$, then $D_h(z,R)$ intersects $Z'$. If $z_k\notin Z'$, then, by the construction of $Z'$, there is some $z_m\in Z'$ so that $z_k\in D_h(z_m,R)$ and hence $z_m\in D_h(z,2R)$. Therefore, in any case, $D_h(z,2R)$ intersects $Z'$ and so $Z'$ is $2R$-dense. If we prove that $Z'$ is sampling for $\U$, then obviously $Z$ is also sampling for $U$. Now, renaming $2R$ as $R$, we have reduced our problem to showing that $Z$ is sampling if it is $R$-dense and the hyperbolic discs with centers from $Z$ and radius $\theta_0R$, where $\theta_0=\frac 14$, are disjoint.
	
\noindent Now let $0<\theta\leq\theta_0$ and let $M$ be a large integer. 
	
\noindent We consider all the points $z_{j,M}$, $1\leq j\leq p(M)$, of $Z$ in $D_h(0,(M+1)R)$. Then for every $z\in D_h(z_{j,M},\theta R)$ we have $f_k(z)\in D_h(f_k(z_{j,M}),\theta N(f_k) R)$. Bearing in mind that for any distinct $z_{j,M},z_{i,M}$ we have 
$$\beta(f(z_{j,M}),f(z_{i,M}))<\frac{1}{k}\,N(f_k)\beta(z_{j,M},z_{i,M})\leq\frac{2(M+1)}{k}\,N(f_k)R,$$ 
we easily see that 
$$ f_k(z)\in D_h\Big(f_k(z_{1,M}),\Big(\theta+\frac{2(M+1)}k\Big) N(f_k)R\Big) $$
for every $z\in\cup_{j=1}^{p(M)}D(z_{j,M},\theta R)$.
	
\noindent On the other hand, for any $z$ on the circle $C_h(0,MR)$ we consider some $z_{j,M}$ with $\beta(z,z_{j,M})<R$ and, in the same manner, we get 
$$f_k(z)\in D_h\Big(f_k(z_{1,M}),\Big(1+\frac{2(M+1)}k\Big) N(f_k)R\Big).$$
These last two statements can be written in the equivalent forms
\begin{equation}\label{ineq1} 
\Big|\frac{f_k(z)-f_k(z_{1,M})}{1-\overline{f_k(z_{1,M})}f_k(z)}\Big|<\frac{e^{(\theta+\frac{2(M+1)}k) N(f_k)R}-1}{e^{(\theta+\frac{2(M+1)}k) N(f_k)R}+1}
\end{equation}
and
\begin{equation}\label{ineq2}
\Big|\frac{f_k(z)-f_k(z_{1,M})}{1-\overline{f_k(z_{1,M})}f_k(z)}\Big|<\frac{e^{(1+\frac{2(M+1)}k) N(f_k)R}-1}{e^{(1+\frac{2(M+1)}k) N(f_k)R}+1}
\end{equation}
respectively.
	
\noindent Now we consider the domain
$$\Omega_M=D_h(0,MR)\setminus \bigcup_{j=p(M_1)+1}^{p(M)}D(z_{j,M},\theta R),$$
where $M_1=\big[\frac{M}{2}\big]+1$ and $z_{j,M}$, $1\leq j\leq p(M_1)$, are all the points of $Z$ which belong to $D_h(0,M_1R)$.
	
\noindent \textit{First claim.} The harmonic measure $\omega(z,E_M,\Omega_M)$, where $E_M$ is the part of the boundary of $\Omega_M$ which is on the circle $C_h(0,MR)$, satisfies
$$\omega(z,E_M,\Omega_M)\leq (1-\omega(\theta,R))^{cM}$$
for every $z\in D_h\big(0,\frac{MR}2\big)$, where $c>0$ is an absolute constant and $0<\omega(\theta,R)<1$.
	
\noindent To prove the claim, we take any $z$ on the circle $C_h(0,(M-6)R)$. Let $z'$ be the closest point of $Z$ to $z$, and $z''$ be the next closest point of $Z$ to $z$. Then $\beta(z',z)<R$ and $\beta(z'',z)<2R$. Since $D_h(z',\theta_0R)$ and $D_h(z'',\theta_0R)$ are disjoint, either $\beta(z',z)\geq\theta_0R$ or $\beta(z'',z)\geq\theta_0R$. So there is $z_1\in Z$ such that $\theta_0R\leq\beta(z_1,z)<2R$. Now it is easy to see that $D_h(z_1,4R)\subseteq D_h(0,MR)$. Indeed, if $w\in D_h(z_1,4R)$ then
$$\beta(w,0)\leq\beta(w,z_1)+\beta(z_1,z)+\beta(z,0)<4R+2R+(M-6)R=MR.$$
Now let $\Omega'=D_h(0,MR)\setminus D_h(z_1,\theta R)$. Then the boundary of $\Omega'$ consists of the circle $E'=C_h(0,MR)$ and of the circle $E''=C_h(z_1,\theta R)$. It is easy to see that
$$\omega(z,E_M,\Omega_M)\leq\omega(z,E',\Omega')=1-\omega(z,E'',\Omega').$$
Moreover
\begin{align*}
\omega(z,E'',\Omega')&\geq\omega(z,C_h(z_1,\theta R),D_h(z_1,4R))\notag\\
&\geq\omega(2R,C_h(0,\theta R),D_h(0,4R))
\end{align*}
by conformal invariance of harmonic measure. Now by the simple formula for harmonic measure for an annulus, we have
$$\omega(2R,C_h(0,\theta R),D_h(0,4R))=\omega(\theta,R),$$
where $\omega(\theta,R)=\frac{\log\big(\frac{e^{4R}-1}{e^{4R}+1}\,\frac{e^{2R}+1}{e^{2R}-1}\big)}{\log\big(\frac{e^{4R}-1}{e^{4R}+1}\,\frac{e^{\theta R}+1}{e^{\theta R}-1}\big)}$.
	
\noindent Thus
$$\omega(z,E_M,\Omega_M)\leq 1-\omega(\theta,R)$$
for any $z$ on the circle $C_h(0,(M-6)R)$. Therefore, by the probabilistic interpretation of harmonic measure for instance, we get
$$\omega(z,E_M,\Omega_M)\leq \omega(z,E_{M-6},\Omega_{M-6})(1-\omega(\theta,R))$$
for $z\in D_h\big(0,\frac{MR}2\big)$.
	
\noindent By induction, for the integer $k$ with $\frac M{12}-1<k\leq\frac M{12}$ we get
$$\omega(z,E_M,\Omega_M)\leq \omega(z,E_{M-6k},\Omega_{M-6k})(1-\omega(\theta,R))^k\leq (1-\omega(\theta,R))^k$$
for $z\in D_h\big(0,\frac{MR}2\big)$, since $D_h\big(0,\frac{MR}2\big)\subseteq\Omega_{M-6k}$, and our first claim is now proved.
	
\noindent Therefore, by (\ref{ineq1}) and (\ref{ineq2}) we get that
$$\Big|\frac{f_k(z)-f_k(z_{1,M})}{1-\overline{f_k(z_{1,M})}f_k(z)}\Big|<A$$
for every $z\in D_h\big(0,\frac{MR}2\big)$, where
\begin{align*}
A=&\frac{e^{(\theta+\frac{2(M+1)}k) N(f_k)R}-1}{e^{(\theta+\frac{2(M+1)}k) N(f_k)R}+1}(1-(1-\omega(\theta,R))^{cM})\\
&+\frac{e^{(1+\frac{2(M+1)}k) N(f_k)R}-1}{e^{(1+\frac{2(M+1)}k) N(f_k)R}+1}(1-\omega(\theta,R))^{cM}.
\end{align*}
This can be written in the form $\beta(f_k(z),f_k(z_{1,M}))<\log\frac{1+A}{1-A}$ and since it holds also for $z=0$, and $f_k(0)=0$, we have that $\beta(f_k(z),0)\leq 2\log\frac{1+A}{1-A}$ and hence
$$|f_k(z)|<\frac{\big(\frac{1+A}{1-A}\big)^2-1}{\big(\frac{1+A}{1-A}\big)^2+1}=\frac{2A}{1+A^2}$$
for every $z\in D_h\big(0,\frac{MR}2\big)$. Since $|z|=\frac{e^{\frac{MR}2}-1}{e^{\frac{MR}2}+1}$ when $\beta(z,0)=\frac{MR}2$, the Lemma of Schwarz implies
\begin{equation}\label{ineq3}
\frac 12 N(f_k)\leq |f_k'(0)|<\frac{e^{\frac{MR}2}+1}{e^{\frac{MR}2}-1}\,\frac{2A}{1+A^2}.
\end{equation}
Now we choose $k=\big[\frac{2(M+1)}{\theta}\big]+1$ and we get
$$A\leq \frac{e^{2\theta N(f_k)R}-1}{e^{2\theta N(f_k)R}+1}(1-(1-\omega(\theta,R))^{cM})+\frac{e^{(1+\theta) N(f_k)R}-1}{e^{(1+\theta) N(f_k)R}+1}(1-\omega(\theta,R))^{cM}.$$

\noindent \textit{Second claim.} We can choose $M$ large enough, independently of $k$, so that
\begin{equation}\label{ineq4}
A\leq\frac{e^{3\theta N(f_k)R}-1}{e^{3\theta N(f_k)R}+1}.
\end{equation}
Comparing the last two inequalities and dividing by $x=N(f_k)R$, we see that it is enough to prove
\begin{align*}
\frac 1x\,\frac{e^{2\theta x}-1}{e^{2\theta x}+1}(1-(1-\omega(\theta,R))^{cM})&+\frac 1x\,\frac{e^{(1+\theta)x}-1}{e^{(1+\theta)x}+1}(1-\omega(\theta,R))^{cM}\notag\\
& \leq\frac 1x\,\frac{e^{3\theta x}-1}{e^{3\theta x}+1}
\end{align*}
when $M$ large enough, independently of $x$, $0<x\leq R$.
	
\noindent Since the function $\frac 1x\,\frac{e^{ax}-1}{e^{ax}+1}$, $a>0$, is decreasing, taking the limit as $x\to 0+$ of the left side of the last inequality and the value of its right side at $x=R$ we see that it is enough to prove
$$\theta(1-(1-\omega(\theta,R))^{cM})+\frac{1+\theta}2(1-\omega(\theta,R))^{cM}\leq\frac 1R\,\frac{e^{3\theta R}-1}{e^{3\theta R}+1}.$$
Hence it is enough to prove
$$\frac{1+\theta}2(1-\omega(\theta,R))^{cM}\leq\frac 1R\,\frac{e^{3\theta R}-1}{e^{3\theta R}+1}-\theta,$$
and so it is enough to prove that
$$\frac 1{\theta R}\,\frac{e^{3\theta R}-1}{e^{3\theta R}+1}>1.$$
Taking the limit of the left side as $\theta R\to 0+$ we find its supremum, and this is $\frac 32$. Therefore, when $\theta\leq\frac{c'}R$, where $c'>0$ is an absolute constant, our last inequality is true and so we proved our second claim.
	
\noindent We are at the point where we know that if $\theta\leq\frac{c'}R$ and $M$ is large enough, depending on $\theta$ and $R$, then (\ref{ineq4}) holds and hence
$$\frac{2A}{A^2+1}\leq \frac{e^{6\theta N(f_k)R}-1}{e^{6\theta N(f_k)R}+1}.$$
Now (\ref{ineq3}) implies
$$\frac 12 N(f_k)<\frac{e^{\frac{MR}2}+1}{e^{\frac{MR}2}-1}\,\frac{e^{6\theta N(f_k)R}-1}{e^{6\theta N(f_k)R}+1},$$
and, taking $M\to\infty$, we find
$$\frac 12 N(f_k)\leq\frac{e^{6\theta N(f_k)R}-1}{e^{6\theta N(f_k)R}+1},$$
when $\theta\leq\frac{c'}R$. Taking again the limit of $\frac 1x\,\frac{e^{6\theta xR}-1}{e^{6\theta xR}+1}$ when $x\to 0+$ we see that
$$\frac 12 \leq \frac 1{N_k(f)}\,\frac{e^{6\theta N(f_k)R}-1}{e^{6\theta N(f_k)R}+1}\leq 3\theta R.$$
When $\theta$ is small enough we get a contradiction.  
\end{proof}

\section{Interpolating sequences of order $n$ for $\U$.}

 The Nevanlinna-Pick interpolation problem says: Given a sequence $\{z_j\}$ of distinct points in $\D$ and another sequence $\{w_j\}$ in $\D$, does there exist $f\in\U$ so that $f(z_j)=w_j$, $j\in\N$?

Nevanlinna in \cite{Nevan} and Pick in \cite{Pick} proved independently that this problem has a solution if and only if for all $N\in\N$ the matrix
$$\left(\frac{1-w_j\overline w_i}{1-z_j\overline z_i}\right)_{j,i=1,\dots,N}$$
is positive semidefinite. This theorem is the root of a very active field connected with many other topics, which may be found for example in \cite{AgMc}. However, the matrix condition is not easy to compute in general, and it does not give information about the geometry of $\{z_j\}$.
 
The Nevanlinna-Pick problem may be seen as a particular case of Carleson's theorem on interpolating sequences for $\Hi$, which we recall now. A sequence $\{z_j\}$ in $\D$ is called an interpolating sequence if for every bounded sequence $\{w_j\}$ there exists $f\in\Hi$ such that $f(z_j)=w_j$, $j\in\N$. In \cite{Carleson} Carleson proved that $\{z_j\}$ is an interpolating sequence if and only if $\{z_j\}$ is separated and there exists $M>0$ such that
$$\sum_{z_j\in Q}(1-|z_j|)\leq M\ell(Q)$$
for any Carleson square $Q$. A Carleson square is a set $Q$ of the form
\begin{equation*}
Q=\left\{r e^{i\theta}:\ 0<1-r<\ell(Q),\,|\theta-\theta_0|<\ell(Q)\right\}.
\end{equation*}
A sequence $\{z_j\}$ in $\D$ is called separated, with constant of separation $\eta$, if 
$$\displaystyle\inf_{i\neq j}\beta(z_i,z_j)=\eta>0.$$ 

This geometric description of interpolating sequences has had a wide impact in Complex Analysis during the last decades, since it has offered a method of studying many different interpolation problems, as one may check for example in \cite{Seip}.

In \cite{MMN} the authors considered a situation which may be understood as intermediate between the Nevanlinna-Pick and the Carleson interpolation problems. They define a sequence $Z=\{z_j\}$ in $\D$ to be an interpolating sequence for $\U$ if there exists $\eps>0$, only depending on $Z$, such that for any sequence $W=\{w_j\}$ in $\D$ satisfying the compatibility condition
\begin{equation}\label{cc1}
\beta(w_i,w_j)\leq \eps \beta(z_i,z_j),\ \ \ i,j\in\N,
\end{equation}
there exists $f\in\U$ such that $f(z_j)=w_j$, $j\in\N$. With this definition the following characterization was proved.

\begin{Talpha}\label{th_MMN}
A sequence $Z$ in $\D$ is an interpolating sequence for $\U$ if and only if the following two conditions hold:
\begin{itemize}  
\item[(a)] $Z=Z^{(1)}\cup Z^{(2)}$, where each $Z^{(i)}$ is a separated sequence.
\item[(b)] There exist $M>0$ and $\alpha$ with $0<\alpha<1$ such that for any Carleson square $Q$ one has
\begin{equation*}
\#\left(Z\cap\left\{z\in Q\, :\, 2^{-m-1}\ell(Q)<1-|z|\leq2^{-m}\ell(Q)\right\}\right)\leq M2^{\alpha m}
\end{equation*}
for any $m=0,1,\dots$.
\end{itemize}
\end{Talpha}

In \cite{MMN} the authors explain that the main condition in the description here is the density condition (b), while the separation 
condition (a) appears because the problem is defined in terms of first differences. Then, one may wonder what may happen to this result if one modifies the definition using higher order differences. To this aim, we will need an appropriate generalization of the Schwarz-Pick Lemma, which is presented below.

In \cite{BM} one may find a version of the classical Schwarz-Pick Lemma involving three points. This was extended in \cite{BRW} for $n$ points, doing a simple iteration of the result in \cite{BM}. Both works pointed out the analogies between the role played by polynomials in the Euclidean setting and finite Blaschke products in the hyperbolic setting. Recall that a finite Blaschke product of degree $n$ is a function of the form
$$\prod_{j=1}^n\frac{|z_j|}{z_j}\frac{z_j-z}{1-\overline z_j z},$$
where $z_j\in\D$, $1\leq j\leq n$. In order to state this result, Beardon and Minda first defined the complex pseudohyperbolic distance between $z,w\in\D$ as follows:
$$[z,w]=\frac{w-z}{1-\overline w z}.$$
Observe that $|[z,w]|=\rho(z,w)$. For a fixed $w$, $[z,w]$ is an automorphism of the disc.

For a fixed $z_1\in\D$ and $f\in\U$ we define the hyperbolic difference quotient as
$$\Delta f(z;z_1)=\begin{cases}
\dfrac{[f(z),f(z_1)]}{[z,z_1]}, &\text{if}\,\, z\in\D\setminus\{z_1\},\\
f^h(z_1), &\text{if}\,\, z=z_1,
\end{cases}$$
where $f^h(z_1)$ is obtained as a limit and represents the hyperbolic derivative of $f$ at $z_1$, that is,
$$f^h(z_1)=\lim_{z\to z_1}\dfrac{[f(z),f(z_1)]}{[z,z_1]}=\frac{(1-|z_1|^2)f'(z_1)}{1-|f(z_1)|^2}.$$
The expression $\Delta f(z;z_1)$ defines a function in $\U$, since it is analytic and, as a consequence of the Swcharz-Pick Lemma, we have $|\Delta f(z;z_1)|\leq1$. 
We may now iterate this process to get differences of higher order. Writing $\Delta^0f(z)=f(z)$, we fix $z_1,\dots,z_n$ in $\D$ and for $k=1,\dots,n$ we define the $k$-th hyperbolic difference quotient as follows:
$$\Delta^k f(z;z_1,\dots,z_k)=\frac{[\Delta^{k-1}f(z;z_1,\dots,z_{k-1}),\Delta^{k-1}f(z_k;z_1,\dots,z_{k-1})]}{[z,z_k]},$$
interpreted as a limit when $z=z_k$. Clearly
$$\Delta^1(\Delta^{k-1}f(\cdot;z_1,\dots,z_{k-1}))(z;z_k)=\Delta^k f(z;z_1,\dots,z_k).$$

The following multi-point Schwarz-Pick Lemma appeared in \cite{BRW}.
\begin{Talpha}\label{SP_general}
Fix $z_1,\dots,z_k$ in $\D$. Then, for all $f\in\U$ and $v,w$ in $\D$,
$$\beta(\Delta^k f(v;z_1,\dots,z_k),\Delta^k f(w;z_1,\dots,z_k))\leq\beta(v,w).$$
Equality holds for a pair (and then for all pairs) of distinct points $v,w$ in $\D$ if and only if $f$ is a Blaschke product of degree $k+1$.
\end{Talpha}

Consider now two sequences $Z=\{z_j\}$ and $W=\{w_j\}$ in $\D$. In the line of \cite{MMN}, we want to define an interpolation problem in $\U$ and give conditions on $Z$ so that there exists a function in $\U$ interpolating $W$ at $Z$. To this end, let us define the hyperbolic difference quotients for the sequence $W$. We fix $z_1,\dots,z_n$ in $Z$ and the corresponding $w_1,\dots,w_n$ in $W$. Writing $\Delta^0_j=w_j$, we define
$$\Delta^k_j=\frac{[\Delta^{k-1}_j,\Delta^{k-1}_k]}{[z_j,z_k]}\,\,\,\text{for}\,\,1\leq k\leq n-1\,\,\text{and}\,\, k+1\leq j\leq n.$$
It is easy to see that each $\Delta^k_j$ depends on $z_1,\dots,z_k,z_j$ and $w_1,\dots,w_k,w_j$. Moreover, if $f$ is a solution of the Nevanlinna-Pick problem then
\begin{equation}\label{estab}
\Delta^k_j=\Delta^kf(z_j;z_1,\dots,z_k)\,\,\,\text{for}\,\,0\leq k\leq n-1\,\,\text{and}\,\, k+1\leq j\leq n,
\end{equation}
(see \cite[Lemma 4.2]{BRW}) so necessarily $|\Delta^k_j|\leq1$ by Theorem \ref{SP_general}.

The complete list of hyperbolic difference quotients may be represented in the following table:
\begin{displaymath}
\begin{array}{cc|llccc|}
\cline{3-7} 
& & \multicolumn{5}{c|}{\text{hyperbolic difference quotients}}\\
\hline
\multicolumn{1}{|c|}{Z} & \multicolumn{1}{|c|}{W} & 1 & 2 & \cdots & n-2 & n-1 \\
\hline 
\multicolumn{1}{|c|}{\vspace{0cm}} & \multicolumn{1}{|c|}{\vspace{0cm}} & & & & & \\
\multicolumn{1}{|c|}{z_1} & \multicolumn{1}{|c|}{w_1=\Delta^0_1} & & & & & \\ 
\multicolumn{1}{|c|}{\vspace{0cm}} & \multicolumn{1}{|c|}{\vspace{0cm}} & & & & & \\
\multicolumn{1}{|c|}{z_2} & \multicolumn{1}{|c|}{w_2=\Delta^0_2} & \Delta^1_2 & & & & \\
\multicolumn{1}{|c|}{\vspace{0cm}} & \multicolumn{1}{|c|}{\vspace{0cm}} & & & & & \\
\multicolumn{1}{|c|}{z_3} & \multicolumn{1}{|c|}{w_3=\Delta^0_3} & \Delta^1_3 & \Delta^2_3 & & &\\
\multicolumn{1}{|c|}{\vspace{0cm}} & \multicolumn{1}{|c|}{\vspace{0cm}} & & & & & \\
\multicolumn{1}{|c|}{\vdots} & \multicolumn{1}{|c|}{\vdots} & \multicolumn{1}{c}{\vdots} & \multicolumn{1}{c}{\vdots} & \ddots & &\\
\multicolumn{1}{|c|}{\vspace{0cm}} & \multicolumn{1}{|c|}{\vspace{0cm}} & & & & & \\
\multicolumn{1}{|c|}{z_{n-1}} & \multicolumn{1}{|c|}{w_{n-1}=\Delta^0_{n-1}} & \Delta^1_{n-1} & \Delta^2_{n-1} & \cdots & \Delta^{n-2}_{n-1} & \\
\multicolumn{1}{|c|}{\vspace{0cm}} & \multicolumn{1}{|c|}{\vspace{0cm}} & & & & & \\
\multicolumn{1}{|c|}{z_{n}} & \multicolumn{1}{|c|}{w_{n}=\Delta^0_{n-1}} & \Delta^1_{n} & \Delta^2_{n} & \cdots & \Delta^{n-2}_{n} & \Delta^{n-1}_n\\
\multicolumn{1}{|c|}{\vspace{0cm}} & \multicolumn{1}{|c|}{\vspace{0cm}} & & & & & \\
\hline
\end{array}
\end{displaymath}
We remark that the hyperbolic difference quotients depend on the order given to the points $z_1,\dots,z_n$. The triangle formed in the table will be called triangle of hyperbolic difference quotients. 

In \cite{BRW} the following generalization of the Nevanlinna-Pick Theorem was proved.

\begin{Talpha}\label{NP_BRW}
Fix distinct $z_1,\dots,z_n$ in $\D$ and corresponding $w_1,\dots,w_n$ in $\D$. Then the Nevanlinna-Pick problem has 
infinitely many solutions if and only if one (and then all) of the following conditions hold:
\begin{enumerate}[(i)]
\item $|\Delta^{n-1}_n|<1$.
\item $|\Delta^{k}_{k+1}|<1$ if $1\leq k\leq n-1$.
\item $|\Delta^{k}_j|<1$ if $1\leq k<j\leq n$.
\item $\beta(\Delta^{k-1}_j,\Delta^{k-1}_k)<\beta(z_j,z_k)$ if $1\leq k<j\leq n$.
\item $\beta(\Delta^{k-1}_j,\Delta^{k-1}_i)<\beta(z_j,z_i)$ if $1\leq k<i<j\leq n$.
\end{enumerate}
\end{Talpha}

We remark that a proof of (ii) in Theorem \ref{NP_BRW} may be essentially found also in \cite[Chapter X]{Wal}. Moreover, this book refers to a paper by Denjoy \cite{Den} where a condition for the infinite points case was proved, using the hyperbolic difference quotients: The Nevanlinna-Pick problem has infinitely many solutions if and only if
$$\sum_{n=1}^\infty\frac{1-|z_n|}{1-|\Delta^{n-1}_n|}<\infty.$$

Fix a sequence $Z$ in $\D$ and $\eps>0$. We will say that a sequence $W$ in $\D$ satisfies the $\eps$-compatibility condition of order $n-1$ for $Z$ if for any $\{z_1,\dots,z_n\}\subset Z$ and the corresponding $\{w_1,\dots,w_n\}\subset W$ the terms of the corresponding triangle of hyperbolic difference quotients satisfy the inequality
\begin{equation}\label{ecc}
 \beta(\Delta^{k}_i,\Delta^{k}_j)\leq\eps\beta(z_i,z_j)\,\,\,\text{for}\,\,0\leq k\leq n-2\,\,\text{and}\,\, k+1\leq i,j\leq n.
\end{equation}
The constant $\eps$ will be called interpolation constant.

Observe that the definition of interpolating sequence given in \cite{MMN} may be rewritten as 
$$\beta(\Delta^0_i,\Delta^0_j)\leq \eps\beta(z_i,z_j),\,\,\, i,j\in\N.$$ 

Now we can give a definition of interpolating sequence, based on the triangle of hyperbolic difference quotients, analogous to the one given in \cite{MMN}.

\begin{defi}\label{def2}
A sequence $Z=\{z_j\}$ in $\D$ is an interpolating sequence of order $n-1$ for $\U$ if there exists $\eps>0$, depending only on the sequence and $n$, such that for any sequence $W=\{w_j\}$ in $\D$ satisfying the $\eps$-compatibility condition \textit{(\ref{ecc})} there exists $f\in\U$ such that $f(z_j)=w_j$, $j\in\N$.
\end{defi}

The definition coincides with the one given in \cite{MMN} when $n=2$.
 
Like the definition in \cite{MMN}, this definition is conformally invariant, in the sense that if $Z=\{z_j\}$ is an interpolating sequence of order $n-1$ and $\tau$ is an automorphism of the disc then $\tau(Z)=\{\tau(z_j)\}$ is also an interpolating sequence of order $n-1$ with the same interpolation constant $\eps$. This is a consequence of the following property which may be found in \cite[Lemma 3.3]{Riv}: Let $f\in\U$, let $\tau$ and $\sigma$ be two automorphisms of the disc and take pairwise distinct points $z_1,\dots,z_n$ in $\D$. Then for $1\leq j\leq n$ there exists $\theta\in [0,2\pi)$ such that 
$$\Delta^j(\sigma\circ f\circ\tau)(z;z_1,\dots,z_j)=e^{i\theta}\Delta^j(f(\tau(z));\tau(z_1),\dots,\tau(z_j)),\,\, z\in\D.$$ 

It is important to remark a direct consequence of definition \ref{def2}: if $Z$ is an interpolating sequence of order $k$ for $\U$, then it is also an interpolating sequence of any order $j\leq k$ for $\U$.

\begin{lma}\label{lemma1}
Let $f$ be an analytic function in $D_h(z_0,\eta_1)\subseteq\D$ such that $|f(z)|\leq C$ on $D_h(z_0,\eta_1)$. Let $a\in D_h(z_0,\eta)$ with fixed $0<\eta<\eta_1$. Then
$$|f(z)-f(a)|\leq \widetilde{C}\rho(z,a)$$
for $z\in D_h(z_0,\eta_1)$, where $\widetilde{C}>0$ depends on $\eta_1-\eta$ and $C$.
\end{lma}

\begin{proof}
A simple consequence of the Schwarz-Pick Lemma.
\end{proof}

\begin{lma}\label{remark}
Let $Z=\{z_1,\dots,z_n\}$ be in a small hyperbolic disc of radius $\eta>0$, and let $W=\{w_1,\dots,w_n\}=\{\Delta^0_1,\dots,\Delta^0_n\}$. Then there exists $C>0$, depending only on $\eta$, such that if for a certain small enough $\eps>0$ one has
\begin{equation}\label{column_cond}
\beta(\Delta^k_j,\Delta^k_{k+1})\leq\eps\beta(z_j,z_{k+1})
\end{equation}
for $1\leq k\leq n-2$ and $k+2\leq j\leq n-1$, then
$$\beta(\Delta^k_i,\Delta^k_j)\leq C\eps\beta(z_i,z_j)$$
for $1\leq k\leq n-2$ and $k+2\leq i<j\leq n-1$.
\end{lma}

\begin{proof}
We shall follow the ideas of the proof of Theorem \ref{NP_BRW}. We consider the following triangle of hyperbolic difference quotients (we shall refer to the last row in a while):
\begin{displaymath}
\begin{array}{cc|lllllc|}
\cline{3-8} 
& & \multicolumn{6}{c|}{\text{hyperbolic difference quotients}}\\
\hline
\multicolumn{1}{|c|}{Z} & \multicolumn{1}{|c|}{W} & 1 & 2 & \cdots & n-2 & n-1 & n \\
\hline 
\multicolumn{1}{|c|}{\vspace{0cm}} & \multicolumn{1}{|c|}{\vspace{0cm}} & & & & & & \\
\multicolumn{1}{|c|}{z_1} & \multicolumn{1}{|c|}{\Delta^0_1} & & & & & & \\ 
\multicolumn{1}{|c|}{\vspace{0cm}} & \multicolumn{1}{|c|}{\vspace{0cm}} & & & & & & \\
\multicolumn{1}{|c|}{z_2} & \multicolumn{1}{|c|}{\Delta^0_2} & \Delta^1_2 & & & & & \\
\multicolumn{1}{|c|}{\vspace{0cm}} & \multicolumn{1}{|c|}{\vspace{0cm}} & & & & & & \\
\multicolumn{1}{|c|}{z_2} & \multicolumn{1}{|c|}{\Delta^0_3} & \Delta^1_3 & \Delta^2_3 & & & & \\
\multicolumn{1}{|c|}{\vspace{0cm}} & \multicolumn{1}{|c|}{\vspace{0cm}} & & & & & & \\
\multicolumn{1}{|c|}{\vdots} & \multicolumn{1}{|c|}{\vdots} & \multicolumn{1}{c}{\vdots} & \multicolumn{1}{c}{\vdots} & \ddots & & & \\
\multicolumn{1}{|c|}{\vspace{0cm}} & \multicolumn{1}{|c|}{\vspace{0cm}} & & & & & & \\
\multicolumn{1}{|c|}{z_{n-1}} & \multicolumn{1}{|c|}{\Delta^0_{n-1}} & \Delta^1_{n-1} & \Delta^2_{n-1} & \cdots & \Delta^{n-2}_{n-1} & & \\
\multicolumn{1}{|c|}{\vspace{0cm}} & \multicolumn{1}{|c|}{\vspace{0cm}} & & & & & & \\
\multicolumn{1}{|c|}{z_{n}} & \multicolumn{1}{|c|}{\Delta^0_n} & \Delta^1_{n} & \Delta^2_{n} & \cdots & \Delta^{n-2}_{n} & \Delta^{n-1}_{n} & \\
\multicolumn{1}{|c|}{\vspace{0cm}} & \multicolumn{1}{|c|}{\vspace{0cm}} & & & & & & \\
\multicolumn{1}{|c|}{z} & \multicolumn{1}{|c|}{g(z)} & \Delta^1g(z) & \Delta^2g(z) & \cdots & \Delta^{n-2}g(z) & \Delta^{n-1}g(z) & \Delta^{n}g(z)\\
\multicolumn{1}{|c|}{\vspace{0cm}} & \multicolumn{1}{|c|}{\vspace{0cm}} & & & & & & \\
\hline
\end{array}
\end{displaymath}
Since $\{z_1,\dots,z_n\}$ is in a small hyperbolic disc $D_h(z,\eta)$, we can consider the pseudohyperbolic distance $\rho$ instead of $\beta$. Hence (\ref{column_cond}) implies
$$|\Delta^k_j|\leq c'\eps\,\,\,\text{for}\,\,1\leq k\leq n-1\,\,\text{and}\,\, k+1\leq j\leq n,$$
where $c'>0$ depends only on $\eta$.
	 
\noindent In particular, if $\epsilon$ is small enough, then all $|\Delta^k_j|$ are smaller than 1, so by Theorem \ref{NP_BRW} there are infinitely many solutions of the Nevanlinna-Pick interpolation problem for $Z$ and $W$. Actually, in the proof appearing in \cite{BRW} we can find the following version of the Schur's algorithm to construct these solutions (see \cite{BRW} for the details). The last row of the table, where 
$$\Delta^kg(z)=\Delta^kg(z;z_1,\dots,z_k)$$ 
for $1\leq k\leq n$, has been added in order to describe this process. 
	
\noindent We take an arbitrary function $g_0\in\U$ and we apply the recursive formula
\begin{equation}\label{formula_NP}
g_k(z)=\left[[z,z_{n+1-k}] g_{k-1}(z),\Delta^{n-k}_{n+1-k}\right]\,\,\,\text{for}\,\,1\leq k\leq n
\end{equation}
to obtain a solution of the problem for $Z$ and $W$. Denoting $g=g_n\in\U$, we easily see that this is a solution, since $g_n(z_j)=\Delta^0_j$ for $1\leq j\leq n$. Moreover, the functions of the algorithm satisfy 
$$g_{n-k}(z)=\Delta^kg(z)\,\,\,\text{for}\,\,0\leq k\leq n.$$
In particular if we choose $g_0=0$ then it is easy to see that there exists $c''>0$ depending only on $c'$ such that $|g_k(z)|\leq c''\eps$ for $1\leq k\leq n-1$. Now property (\ref{estab}) and Lemma \ref{lemma1} applied in the whole disc $\D$ show that
$$|\Delta^k_j-\Delta^k_i|=|g_{n-k}(z_j)-g_{n-k}(z_i)|\leq C\eps\rho(z_j,z_i)$$
for  $1\leq k\leq n-2$ and $k+1\leq i<j\leq n$, where $C>0$ depends on $c''$ and $\eta$ and hence on $\eta$.
\end{proof}

We observe that Lemma \ref{remark} is stated for a fixed ordering of the $n$ points and for $k\geq 1$. The case $k=0$ holds if we suppose that the $\eps$-compatibility condition is true for any permutation of the points. Hence, in order to check inequality (\ref{ecc}) for $0\leq k\leq n-2$ when $\{z_1,\ldots,z_n\}$ is in a small hyperbolic disc we just need to see that 
\begin{equation}\label{ecc2}
|\Delta^k_j|\leq\eps\,\,\,\text{for}\,\,1\leq k\leq n-1\,\,\text{and}\,\, k+1\leq j\leq n, 
\end{equation}
and we need to check it for any permutation of $\{z_1,\ldots,z_n\}$.

The second result of this work is the characterization of these interpolating sequences.

\begin{thm}\label{main_th}
A sequence $Z$ in $\D$ is an interpolating sequence of order $n-1$ for $\U$ if and only if the following two conditions hold:
\begin{itemize}
\item[(a)] $Z=Z^{(1)}\cup\dots\cup Z^{(n)}$, where each $Z^{(i)}$ is a separated sequence.
\item[(b)] There exist $M>0$ and $\alpha$ with $0<\alpha<1$ such that for any Carleson square $Q$ one has
\begin{equation*}
\#\left(Z\cap\left\{z\in Q\, :\, 2^{-m-1}\ell(Q)<1-|z|\leq2^{-m}\ell(Q)\right\}\right)\leq M2^{\alpha m}
\end{equation*}
for any $m=0,1,\dots$.
\end{itemize}
\end{thm}

\begin{proof}

\noindent \textit{Necessity.}

\noindent Let $Z\subset\D$ be an interpolating sequence of order $n-1$ for $\U$. We need to show that there exists $\eta>0$ such that there cannot be more than $n$ points inside any hyperbolic disc $D_h(z,\eta)$. To this end, take $z_1,z_2,\dots,z_n,z_{n+1}\in Z$ such that $z_j\in\D_h(z_1,\eta)$ for $2\leq j\leq n+1$. Suppose also that 
$$\rho(z_1,z_j)\geq\rho(z_1,z_{n+1})\,\,\,\text{for}\,\,2\leq j\leq n+1.$$

\noindent Now we consider the values $w_j=0$ for $1\leq j\leq n$ and $w_{n+1}=\eps x$, where $x$ will be determined so that the sequence $W=\{w_1,\dots,w_{n+1}\}$ satisfies the $\eps$-compatibility condition. To this end, we have to take the set $Z^j=\{z_1,\dots,z_{n+1}\}\setminus\{z_j\}$ for $1\leq j\leq n+1$ and check inequality (\ref{ecc2}) for the terms in the corresponding triangle of hyperbolic difference quotients for any permutation of $Z^j$. We will see that it is enough to check it only for a certain permutation, with a suitable choice of the value $x$.

\noindent The inequality trivially holds for the set $Z^{n+1}$ since all the hyperbolic differences vanish. For the sets $Z^j$ with $1\leq j\leq n$ we will first consider the original ordering of the points, and then we will consider any permutation of them. For the sake of simplicity, we will do it only in the case of $Z^1$ but the argument is the same for all other sets. Consider the following triangle of hyperbolic difference quotients.
\begin{displaymath}
\begin{array}{cc|lllll|}
\cline{3-7} 
& & \multicolumn{5}{c|}{\text{hyperbolic difference quotients}}\\
\hline
\multicolumn{1}{|c|}{Z} & \multicolumn{1}{|c|}{W} & 1 & 2 & \cdots & n-2 & n-1 \\
\hline 
\multicolumn{1}{|c|}{\vspace{0cm}} & \multicolumn{1}{|c|}{\vspace{0cm}} & & & & & \\
\multicolumn{1}{|c|}{z_2} & \multicolumn{1}{|c|}{0} & & & & & \\ 
\multicolumn{1}{|c|}{\vspace{0cm}} & \multicolumn{1}{|c|}{\vspace{0cm}} & & & & & \\
\multicolumn{1}{|c|}{z_3} & \multicolumn{1}{|c|}{0} & \Delta^1_3 & & & & \\
\multicolumn{1}{|c|}{\vspace{0cm}} & \multicolumn{1}{|c|}{\vspace{0cm}} & & & & & \\
\multicolumn{1}{|c|}{z_4} & \multicolumn{1}{|c|}{0} & \Delta^1_4 & \Delta^2_4 & & & \\
\multicolumn{1}{|c|}{\vspace{0cm}} & \multicolumn{1}{|c|}{\vspace{0cm}} & & & & & \\
\multicolumn{1}{|c|}{\vdots} & \multicolumn{1}{|c|}{\vdots} & \multicolumn{1}{c}{\vdots} & \multicolumn{1}{c}{\vdots} & \ddots & & \\
\multicolumn{1}{|c|}{\vspace{0cm}} & \multicolumn{1}{|c|}{\vspace{0cm}} & & & & & \\
\multicolumn{1}{|c|}{z_{n}} & \multicolumn{1}{|c|}{0} & \Delta^1_{n} & \Delta^2_{n} & \cdots & \Delta^{n-2}_{n} & \\
\multicolumn{1}{|c|}{\vspace{0cm}} & \multicolumn{1}{|c|}{\vspace{0cm}} & & & & & \\
\multicolumn{1}{|c|}{z_{n+1}} & \multicolumn{1}{|c|}{\eps x} & \Delta^1_{n+1} & \Delta^2_{n+1} & \cdots & \Delta^{n-2}_{n+1} & \Delta^{n-1}_{n+1} \\
\multicolumn{1}{|c|}{\vspace{0cm}} & \multicolumn{1}{|c|}{\vspace{0cm}} & & & & & \\
\hline
\end{array}
\end{displaymath}
It is easy to see that the only hyperbolic difference quotients that do not necessarily vanish are the ones in the row corresponding to the point $z_{n+1}$, that is, the ones of the form $\Delta^k_{n+1}$ for $1\leq k\leq n-1$, and moreover,
$$\Delta^k_{n+1}=\frac{\Delta^{k-1}_{n+1}}{[z_{n+1},z_{k+1}]},$$
where $\Delta^0_{n+1}=\eps x$. Consequently, (\ref{ecc2}) is true in this case if
$$|x|\leq\prod_{i=2}^{n}\rho(z_i,z_{n+1}),$$
Without loss of generality, we may suppose that
$$\rho(z_j,z_{n+1})\leq\rho(z_n,z_{n+1})\,\,\,\text{for}\,\,1\leq j\leq n.$$
Then we take 
$$x=C\prod_{j=1}^{n-1}[z_{n+1},z_j]$$ 
with $0<C<1$, and we have that $|\Delta^k_j|\leq C\eps$ for $1\leq k\leq n-1$ and $k+2\leq j\leq n+1$.

\noindent Next, we are going to see that, for a suitable choice of $C$, inequality (\ref{ecc2}) holds for the terms of the triangle corresponding to any permutation of the points. To this end, we will proceed as in the proof of Lemma \ref{remark}. We take $g_0=0$ and, since all the terms in the upper diagonal but the last one vanish, the recursive formula (\ref{formula_NP}) produces the following solution of the interpolation problem for $Z^{1}$:
$$g(z)=(-1)^{n+1}\Delta^{n-1}_{n+1}\prod_{j=2}^n[z,z_{j}].$$
We easily see that $|g(z)|\leq\eps$. Hence, choosing $C$ small enough, Lemma \ref{lemma1} and property (\ref{estab}) show that for any permutation of the points inequality (\ref{ecc2}) holds. The constant $C$ does not depend on the points, so the choice of $x$ is the same for any $Z^j$, $1\leq j\leq n$. 

\noindent Since $Z$ is an interpolating sequence of order $n-1$ there is a function $f\in\U$ such that $f(z_j)=w_j$ for $1\leq j\leq n+1$. Besides, the function has the form
$$f(z)=\prod_{j=1}^{n}[z,z_j]h(z),$$
where $h\in\U$. Now evaluating the function at $z_{n+1}$ we see that
$$\eps\prod_{j=1}^{n-1}\rho(z_j,z_{n+1})=|f(z_{n+1})|\leq\prod_{j=1}^n\rho(z_j,z_{n+1}).$$
Hence, there exists an absolute constant $c>0$ such that $\eta>c\eps$. Therefore, if $\eta\leq c\eps$, we get at a contradiction and we have proved the necessity of (a).

\noindent The necessity of (b) is a consequence of Theorem \ref{th_MMN} since, in particular, each separated sequence $Z^{(j)}$ is an 
interpolating sequence of order 1.

\noindent \textit{Sufficiency.}

\noindent Assume that the sequence $Z$ of points in $\D$ satisfies (a) and (b). We consider $\eps>0$, which will be chosen appropriately later, and we assume that the sequence $W$ in $\D$ satisfies the $\eps$-compatibility condition.

\noindent We take first the separated sequence $Z^{(1)}$ and the corresponding $W^{(1)}$. Condition (b) and the $\eps$-compatibility condition for first differences allows us to construct a function $f_1\in\U$ that interpolates $W^{(1)}$ at $Z^{(1)}$. Furthermore, this function has two additional properties. In order to state them, let $E_1(z)$ be the outer function with boundary values $1-|f_1(e^{i\theta})|$, that is,
$$E_1(z)=\exp\left(\frac{1}{2\pi}\int_0^{2\pi}\frac{e^{i\theta}+z}{e^{i\theta}-z}\log{\left(1-|f_1(e^{i\theta})|\right)}d\theta\right).$$
The construction of $f_1$ is explained in detail in \cite[Section 4]{MMN} and its properties are stated in the following Lemma. 

\begin{Lalpha}
Let $Z^{(1)}=\{z^{(1)}_i\}$ be a separated sequence in $\D$ with interpolation constant $\eps>0$. Let $W^{(1)}=\{w^{(1)}_i\}$ be a sequence in $\D$ such that the compatibility condition (\ref{cc1}) holds. Then there exists $f_1\in\U$ such that $f_1(z^{(1)}_i)=w^{(1)}_i$, $i\in\N$. Moreover, there exists $C>0$ such that 
\begin{equation}\label{prop1}
|E_1(z)|\geq C(1-|f_1(z)|).
\end{equation}
for every $z\in\D$. Also, there exists $\eta_1>0$ depending on $Z^{(1)}$ such that
\begin{equation}\label{prop2}
\beta(f_1(z),f_1(z^{(1)}_i))\leq C\eps\beta(z,z^{(1)}_i)
\end{equation}
for $z\in D_h(z^{(1)}_i,\eta_1)$ and $i\in\N$.
\end{Lalpha}

\noindent Now let $Z^{(j)}=\{z^{(j)}_i\}$, $2\leq j\leq n$. As in \cite{MMN}, we can take $\eta>0$ to be smaller than the separation constant of the sequence $Z^{(1)}$ and also assume that $\eta<\eta_1$, where $\eta_1$ is the radius appearing in (\ref{prop2}). Then we can suppose that 
$$Z^{(j)}\subseteq\bigcup_{i\in\N} D_h(z^{(1)}_i,\eta)\,\,\,\text{for}\,\,2\leq j\leq n.$$
This means that for each $z^{(j)}_m\in Z^{(j)}$ there exists $z^{(1)}_{i(m)}\in Z^{(1)}$ such that $z^{(j)}_m\in D_h(z^{(1)}_{i(m)},\eta)$.

\noindent Now we will state and solve an auxiliary interpolation problem for $Z\setminus Z^{(1)}$ as follows. 

\noindent Let $B_1(z)$ be the Blaschke product with zeros $Z^{(1)}$ and, for $2\leq j\leq n$, let
$$\widetilde{w}^{(j)}_m=\frac{f_1(z^{(j)}_m)-w^{(j)}_m}{B_1(z^{(j)}_m) E_1(z^{(j)}_m)}.$$
These auxiliary values were also defined in \cite{MMN}. Hence we can address the Nevanlinna-Pick interpolation problem for 
$Z\setminus Z^{(1)}$ and $\widetilde{W}=\cup_{j=2}^n\widetilde{W}^{(j)}$, with $\widetilde{W}^{(j)}=\{\widetilde{w}^{(j)}_m\}$. If there exists a solution $\widetilde{f}\in\U$ such that $\widetilde{f}(z^{(j)}_m)=\widetilde{w}^{(j)}_m$ for $2\leq j\leq n$ and $m\in\N$, then the function
\begin{equation}\label{solution}
f(z)=f_1(z)-B_1(z) E_1(z) \widetilde{f}(z)
\end{equation}
is in $\U$ and solves the Nevanlinna-Pick problem for $Z$ and $W$. 

\noindent We will proceed by induction in $n$. The case $n=2$ is actually the proof of the sufficiency of Theorem \ref{th_MMN}. So let us suppose that Theorem \ref{main_th} holds for $n-1$ and prove it for $n$.

\noindent In order to make the proof easier, we will split the auxiliary problem into two different interpolation problems, since
\begin{equation}\label{splitval}
\widetilde{w}^{(j)}_m=\frac{f_1(z^{(j)}_m)-f_1(z^{(1)}_{i(m)})}{B_1(z^{(j)}_m) E_1(z^{(j)}_m)}+\frac{w^{(1)}_{i(m)}-w^{(j)}_m}{B_1(z^{(j)}_m) E_1(z^{(j)}_m)}.
\end{equation}
Consider first the interpolation problem defined by the values
\begin{equation*}
\hdel^{(j)}_m=2\,\frac{f_1(z^{(j)}_m)-f_1(z^{(1)}_{i(m)})}{B_1(z^{(j)}_m) E_1(z^{(j)}_m)},
\end{equation*}
corresponding to $z^{(j)}_m$ for $2\leq j\leq n$ and $m\in\N$. Observe that properties (\ref{prop1}) and (\ref{prop2}) imply that 
$|\hdel^{(j)}_m|\lesssim\eps$. Now we want to apply the hypothesis of induction to get the existence of the function $h\in\U$ such that $h(z^{(j)}_m)=\hdel^{(j)}_m$ for $2\leq j\leq n$ and $m\in\N$. In order to do this, we will choose $n-1$ points of the sequence $Z\setminus Z^{(1)}$ and check the $\eps$-compatibility condition for the corresponding triangle of hyperbolic difference quotients. 

\noindent To this aim we have basically two different cases: when the points are far and when they are close, that is, when the distance between the points is bigger than a fixed constant and when it is smaller. Let us start by the first case, and in order to simplify the notation, take first $z_2,\dots,z_n\in Z\setminus Z^{(1)}$ such that $\beta(z_i,z_j)\geq\eta$ for $2\leq i<j\leq n$. Denote the corresponding images by $\hdel^0_j$ for $2\leq j\leq n$. Then the triangle of hyperbolic difference quotients is 
\begin{displaymath}
\begin{array}{|c|c|cccc|}
\hline
Z & W & 1 & 2 & \cdots & n-2 \\
\hline
\vspace{0cm} & & & & & \\
z_2 & \hdel^0_2 & & & & \\
\vspace{0cm} & & & & & \\
z_3 & \hdel^0_3 & \hdel^1_3 & & & \\
\vspace{0cm} & & & & & \\
z_4 & \hdel^0_4 & \hdel^1_4 & \hdel^2_4 & & \\
\vspace{0cm} & & & & & \\
\vdots & \vdots & \vdots & \vdots & \ddots & \\
\vspace{0cm} & & & & & \\
z_n & \hdel^0_n & \hdel^1_n & \hdel^2_n & \cdots & \hdel^{n-2}_n \\
\vspace{0cm} & & & & & \\
\hline
\end{array}
\end{displaymath}
where
$$\hdel^k_j=\frac{\left[\hdel^{k-1}_j,\hdel^{k-1}_{k+1}\right]}{[z_j,z_{k+1}]}\,\,\,\text{for}\,\,1\leq k\leq n-2\,\,\text{and}\,\,k+2\leq j\leq n.$$
Since $|\hdel^0_j|\lesssim\eps$ if $\eps$ is small enough, one may see that $\beta(\hdel^k_j,\hdel^k_i)\leq C\eps$, so then there exists $C'=C'(\eta)>0$ such that 
$$\beta(\hdel^k_j,\hdel^k_i)\leq C'\eps\beta(z_j,z_i)\,\,\,\text{for}\,\,0\leq k\leq n-2\,\,\text{and}\,\,k+2\leq i<j\leq n.$$
Hence we may focus on the case of close points. We take $z_1\in Z^{(1)}$ so that for each $j=2,\dots,n$ there exists a unique $z_j\in Z^{(j)}\cap D_h(z_1,\eta)$. We consider $\Omega_1=D_h(z_1,\eta_1)$ and, since $\eta_1>\eta$, all $z_j$ are farther than a fixed distance from $\partial\Omega_1=C_h(z_1,\eta_1)$. Then the corresponding images are
$$\hdel^0_j=2\,\frac{f_1(z_j)-f_1(z_1)}{B_1(z_j)E_1(z_j)}\,\,\,\text{for}\,\,2\leq j\leq n.$$
We want to see that $|\hdel^k_j|\leq C \eps$, with $C$ depending only on the sequence $\{z_2,\dots,z_n\}$. We consider the function
$$f_2(z)=2\frac{f_1(z)-f_1(z_1)}{B_1(z)E_1(z)},$$
which by (\ref{prop1}) and (\ref{prop2}) is analytic and $|f_2(z)|\lesssim\eps$ in $\Omega_1$. Moreover, $f_2(z_j)=\hdel^0_j$ for $2\leq j\leq n$. We observe that Lemma \ref{lemma1} implies that $|f_2(z_2)-f_2(z)|\lesssim\eps\rho(z,z_2)$ for $z\in\Omega_1$. Now, the function
$$f_3(z)=\frac{[f_2(z),f_2(z_2)]}{[z,z_2]}$$
is analytic in $\Omega_1$ and $|f_3(z)|\lesssim\eps$. Moreover, $f_3(z_j)=\hdel^1_j$ for $3\leq j\leq n$. Inductively, we consider the function
$$f_{k+2}(z)=\frac{[f_{k+1}(z),f_{k+1}(z_{k+1})]}{[z,z_{k+1}]},$$
which is analytic in $\Omega_1$, bounded by a constant comparable to $\eps$ and $f_{k+2}(z_j)=\hdel^{k}_j$ for $k+2\leq j\leq n$. Then we may apply Lemma \ref{lemma1} to conclude that
$$|\hdel^k_{k+2}-\hdel^k_j|\lesssim\eps\rho(z_{k+2},z_j)\,\,\,\text{for}\,\,k+3\leq j\leq n.$$
Hence, choosing $\eps$ small enough, (\ref{ecc2}) holds for $n-1$ points. Clearly it also holds for any permutation of them. We can consequently apply Theorem \ref{main_th} to assert that there exists $h\in\U$ such that $h(z^{(j)}_m)=\hdel^{(j)}_m$ for $2\leq j\leq n$ and $m\in\N$.

\noindent Now we consider the second interpolation problem defined by the values
\begin{equation*}
\tdel^{(j)}_m=2\,\frac{w^{(j)}_m-w^{(1)}_{i(m)}}{B_1(z^{(j)}_m) E_1(z^{(j)}_m)}.
\end{equation*}
We observe that (\ref{cc1}) and (\ref{prop1}) imply that $|\tdel^{(j)}_m|\lesssim\eps$. As in the previous problem, we can consider two different cases. If $\beta(z^{(j)}_m,z^{(i)}_l)\geq\eta$, arguing as above we can see that for $\eps$ small enough condition (\ref{ecc}) holds. Hence, we may focus again on the domain $\Omega_1$. We observe the following tables of hyperbolic difference quotients
\begin{multicols}{2}
\begin{displaymath}
\begin{array}{|c|c|cccc|}
\hline
\vspace{0cm} & & & & &\\
z_1 & \Delta^0_1 & & & &\\
\vspace{0cm} & & & & &\\
z_2 & \Delta^0_2 & \Delta^1_2 & & & \\
\vspace{0cm} & & & & & \\
z_3 & \Delta^0_3 & \Delta^1_3 & \Delta^2_3 & & \\
\vspace{0cm} & & & & & \\
\vdots & \vdots & \vdots & \vdots & \ddots & \\
\vspace{0cm} & & & & & \\
z_n & \Delta^0_n & \Delta^1_n & \Delta^2_n & \cdots & \Delta^{n-1}_n \\
\vspace{0cm} & & & & & \\
\hline
\end{array}
\end{displaymath} 

\begin{displaymath}
\begin{array}{|c|c|ccc|}
\hline
\vspace{0cm} & & & & \\
z_2 & \tdel^0_2 & & & \\
\vspace{0cm} & & & & \\
z_3 & \tdel^0_3 & \tdel^1_3 & &  \\
\vspace{0cm} & & & & \\
\vdots & \vdots & \vdots & \ddots & \\
\vspace{0cm} & & & & \\
z_n & \tdel^0_n & \tdel^1_n & \cdots & \tdel^{n-2}_n \\
\vspace{0cm} & & & & \\
\hline
\end{array}
\end{displaymath} 
\end{multicols}
where
$$\tdel^0_j=2\,\frac{\Delta^0_1-\Delta^0_j}{B_1(z_1) E_1(z_1)}\,\,\,\text{for}\,\,2\leq j\leq n$$
and
$$\tdel^k_j=\frac{[\tdel^{k-1}_j,\tdel^{k-1}_{k+1}]}{[z_j,z_{k+1}]}\,\,\,\text{for}\,\,1\leq k\leq n-2\,\,\text{and}\,\,k+2\leq j\leq n.$$
We recall that, by hypothesis, $|\Delta^k_j|\leq\eps$ for $1\leq k\leq n-1$ and $k+1\leq j\leq n$. This fact will imply that there exists $C>0$, depending only on $z_1,\dots,z_n$, such that $|\tdel^k_j|\leq C\eps$ for $1\leq k\leq n-2$ and $k+2\leq j\leq n$. To this end, Theorem \ref{NP_BRW} applied to the triangle on the left says that there exists $g\in\U$ such that $g(z_j)=\Delta^0_j$ for $1\leq j\leq n$. Arguing as in the proof of Lemma \ref{remark}, we can take $g_0=0$ and then formula (\ref{formula_NP}) generates a solution $g\in\U$ such that $g(z_j)=w_j$ for $1\leq j\leq n$ and such that $|\Delta^1g(z;z_1)|\leq C\eps$ for $z\in\D$, where $C>0$ is an absolute constant (depending only on $n$).

\noindent Let $B_1^*(z)=\frac{B_1(z)}{[z,z_1]}$. Then there exists $C>0$, depending on $\eta$, such that $|B_1^*(z)|\geq C$ when $z\in \Omega_1$. Besides, we observe that
$$\tdel^0_j=2\,\frac{\Delta^0_1-\Delta^0_j}{B_1(z_j)E_1(z_j)}=2\,\frac{1-\overline{\Delta^0_1}\Delta^0_j}{B_1^*(z_j)E_1(z_j)}\,\Delta^1_j$$
for $2\leq j\leq n$. Then we can define the function
\begin{equation*}
\widetilde{g}(z)=2\,\frac{1-\overline{\Delta^0_1}\,g(z)}{B_1^*(z)E_1(z)}\,\Delta^1g(z;z_1),
\end{equation*}
which is analytic on $\Omega_1$. Since $|E_1(z)|\simeq|1-\overline{\Delta^0_1}g(z)|$ when $z\in \Omega_1$, there exists $C>0$ such that $|\widetilde{g}(z)|\leq C\eps$ on $\Omega_1$. Furthermore, by (\ref{estab}), and choosing $\eps$ small enough, the function $\widetilde{g}$ is in $\U$ and interpolates the values $\tdel^0_j$ at $z_j$ for $2\leq j\leq n$. Now, an argument based on Lemma \ref{lemma1}, similar to the one used in the previous interpolation problem with the functions $f_k$, shows that there exists a constant $\widetilde{C}>0$ depending only on the sequence such that
$$|\tdel^k_j|=|\Delta^k\widetilde{g}(z_j;z_2,\dots,z_{k+1})|\leq\widetilde{C}\eps$$
for $1\leq k\leq n-2$ and $k+2\leq j\leq n$. 

\noindent Consequently, for $\eps>0$ small enough the values $\{\tdel^{(j)}_m\}$ satisfy condition (\ref{ecc2}). And it is clear that this holds for any permutation of $\{z_2,\dots,z_n\}$. Applying now, by the induction hypothesis, Theorem \ref{main_th} to $Z^{(2)}\cup\dots\cup Z^{(n)}$, we get that there exists a function $\widetilde{h}\in\U$ such that $\widetilde{h}(z^{(j)}_m)=\tdel^{(j)}_m$ for $2\leq j\leq n$ and $m\in\N$.

\noindent Finally, by (\ref{splitval}), the function $\widetilde{f}=\dfrac{1}{2}(h+\widetilde{h})$ is in $\U$ and $\widetilde{f}(z^{(j)}_m)=\widetilde{w}^{(j)}_m$ for $2\leq j\leq n$ and $m\in\N$. Hence (\ref{solution}) gives the solution of the problem for $n$ separated sequences, and Theorem \ref{main_th} is proved.
\end{proof}

Observe that condition (b) in Theorem \ref{main_th}, that is, the density condition, remains the same as in Theorem \ref{th_MMN}, hence the 
only change is on the separation condition (a), as it was expected. Furthermore, other previous results on generalizations of interpolation problems for Hardy spaces pointed out to this fact. In \cite{Vas2} and \cite{Vas1} the author generalizes the theorem of Carleson on interpolating sequences for $\Hi$. A suitable notion of hyperbolic difference quotient is given there. Then, imposing the boundedness of the difference quotient up to a certain order, leads to changing the separation condition in Carleson's result. Also in \cite{BNO} the authors generalize the interpolation problem for $\Hp$ solved previously in \cite{SS} (case $p\geq1$) and \cite{Kab} (case $0<p<1$). 
In that work, the trace space is defined by means of a certain maximal function using the hyperbolic difference quotients, and the separation condition changes in the same direction as in our result. These problems are also studied in \cite{Har1} and \cite{Har2}.

\bibliographystyle{plain}
\bibliography{monreal_papadimitrakis}

\end{document}